\documentclass[a4paper,12pt]{article}

\usepackage{amssymb}
\usepackage{graphicx}
\usepackage{amsmath}
\usepackage{amsfonts}

\newcommand{\singlespacing}{\let\CS=\@currsize\renewcommand{\baselinestreatch}{1.0}\tiny\CS}
\newcommand{\doublespacing}{\let\CS=\@currsize\renewcommand{\baselinestreatch}{1.5}\tiny\CS }

\newtheorem{thm}{Theorem}[section]

\newtheorem{definition}{Definition}

\newtheorem{lem}{Lemma}[section]

\numberwithin{equation}{section}

\begin{document}
\begin{center}
\textbf{\Large {Certain Almost Kenmotsu Metrics Satisfying the Miao-Tam Equation}}
\end{center}
\centerline{Dhriti Sundar Patra$^1$ and Amalendu Ghosh$^2$}

\newtheorem{Theorem}{\quad Theorem}[section]
\newtheorem{Definition}[Theorem]{\quad Definition}
\newtheorem{Corollary}[Theorem]{\quad Corollary}
\newtheorem{Lemma}[Theorem]{\quad Lemma}
\newtheorem{Example}[Theorem]{\emph{Example}}
\newtheorem{Proposition}[Theorem]{Proposition}
\numberwithin{equation}{section}
\noindent\\
\textbf{Abstract:} {In this paper we characterize certain class of almost Kenmotsu metrics satisfying the Miao-Tam equation.}



\section{Introduction}
A classical problem in differential geometry is to find Riemannian metrics on a given compact manifold $M^n$ that provides constant scalar curvature. In this sense, it is crucial to study the critical points of the total scalar curvature functional through variational approach. Einstein and Hilbert proved that the critical points of the total scalar curvature functional $\mathcal{S} :\mathcal{M} \longrightarrow  R$ defiend by
\begin{eqnarray*}
\mathcal{S}(g) =\int_{M}r_{g}dv_{g},
\end{eqnarray*}
on a compact orientable Riemannian manifold $(M^n,g)$ (where $\mathcal{M}$ denotes the set of all Riemannian metrics on $(M^n,g)$ of unit volume, $r_{g}$ the scalar curvature and $dv_{g}$ the volume form of $g$) restricted to the set of all Riemannian metrics $\mathcal{M}$ of unity volume are Einstein (see \cite{AB}). This stimulated many interesting research. In \cite{MT} the authors studied the variational properties of the volume functional over the space of constant scalar curvature on a given compact Riemannian manifold with boundary. This leads to the following definition
\begin{definition}
Let $(M^n, g)$, $n>2$ be a compact Riemannian manifold with a smooth boundary metric $\partial M$. Then $g$ is said to be a critical metric if there exists a smooth function $\lambda : M^n \rightarrow \mathbb{R}$ such that
\begin{eqnarray}\label{1.1}
-(\Delta_g\lambda)g + \nabla^{2}_g \lambda - \lambda S = g
\end{eqnarray}
on $M$ and $\lambda = 0$ on $\partial M$, where $ \Delta_g$, $\nabla^{2}_g \lambda$ are the Laplacian, Hessian operator with respect to the metric $g$ and $S$ is the $(0, 2)$ Ricci curvature of $g$. The function $\lambda $ is known as the potential function.
\end{definition}

For brevity, the metrics which satisfies (\ref{1.1}) are known as Miao-Tam critical metrics and we refer the equation (\ref{1.1}) as Miao-Tam equation. In  \cite{MT}, Miao-Tam prove that any Riemannian metric $g$ satisfying the equation (\ref{1.1}) must have constant scalar curvature. From which it follows that a critical metric $g$ always has constant scalar curvature. The existence of such metrics was proved on some certain classes of warped product spaces which include the usual spatial Schwarzschild metrics and Ads-Schwarzschild metrics restricted to certain domains containing their horizon and bounded by two spherically symmetric spheres (cf. Corollary $3.1$ and Corollary $3.2$ of \cite{MTE}).\\

In \cite{MT}, the authors classified Einstein and conformally flat Riemannian manifold satisfying (\ref{1.1}). In fact, they proved that any connected, compact, Einstein manifold with smooth boundary satisfying Miao-Tam critical condition is isometric to a geodesic ball in a simply connected space form $\mathbb{R}^{n}$, $\mathbb{H}^{n}$ or $\mathbb{S}^{n}$. Similar characterization was obtained when $g$ is a conformally flat metric on a simply connected manifold $M$ such that the boundary of $(M, g)$ is isometric to a round sphere. We also note that the last result has been generalized in dimension $4$ under the Bach flat assumption by Barros-Ribeiro \cite{BDR}. Recently, the authors studied equation (\ref{1.1}) on certain class of odd dimensional Riemannian manifold, namely contact metric manifold (see \cite{PG}). Particularly, it was proved that a complete $K$-contact metric satisfying the Miao-Tam critical condition is isometric to a unit sphere $S^{2n+1}$. This result intrigues us to study the Miao-Tam equation on other almost contact metric manifolds. In this paper, we classify certain class of almost Kenmotsu manifold which satisfies Miao-Tam equation.

\section{Preliminaries}
A contact manifold is a Riemannian manifold $M$ of dimension $(2n + 1)$ which carries a global $1$-form $\eta$ such that $\eta\wedge(d\eta)^{n}\ne0$  everywhere on $M$. The form $\eta$ is usually known as the contact form on $M$. It is well known that a contact manifold admits an almost contact metric structure on $(\varphi,\xi,\eta, g)$, where $\varphi$ is a tensor field of type $(1,1)$, $\xi$ a global vector field known as the characteristic vector field (or the Reeb vector field) and $g$ is Riemannian metric, such that
\begin{eqnarray}\label{2.1}
\varphi^2X = -X + \eta(X) \xi, \hskip 0.3cm \eta(\xi)=1,
\end{eqnarray}
\begin{eqnarray}\label{2.2}
g(\varphi X, \varphi Y) = g(X,Y) - \eta(X)\eta(Y),
\end{eqnarray}
for all vector fields $X$, $Y$ on $M$. It follows from equation (\ref{2.1}) that $\varphi\xi=0$ and $\eta \circ\varphi=0$ (see \cite{Blair}, p.43). A Riemannian manifold $M$ together with the almost contact metric structure $(\varphi,\xi,\eta,g)$ is said to be a almost contact metric manifold.
On almost contact metric manifolds one can always define a fundamental $2$-form $\Phi$ by $\Phi(X,Y) = g(X,\varphi Y)$ for all vector fields $X$, $Y$ on $M$. An almost contact metric structure of $M$ is said to be contact metric if $\Phi = d\eta$, and is said to be almost Kenmotsu manifold if $d\eta = 0$ and $d\Phi = 2\eta\wedge\Phi$. Further, a condition for an almost contact metric structure being normal is equivalent to vanishing of the $(1,2)$-type torsion tensor $N_\varphi$, defined by $N_\varphi = [\varphi,\varphi] + 2d\eta \times \xi$, where $[\varphi,\varphi]$ is the Nijenhuis torsion of $\varphi$. A normal almost Kenmotsu manifold is said to be a Kenmotsu manifold and the normality condition is given by
$(\nabla_{X}\varphi)Y = g(\varphi X,Y)\xi - \eta(Y)\varphi X$
for all vector fields $X$, $Y$ on $M$. In \cite{K}, Kenmotsu proved that a warped product of a line and a K\"{a}hlerian manifold admits a Kenmotsu structure. In fact, a Kenmotsu manifold $M$ is locally a warped product $I\times_{f}M^{2n}$, where $I$ is an open interval with coordinate $t$, $f = ce^{2t}$ is the warping function for some positive constant $c$ and $M^{2n}$ is a K\"{a}hlerian manifold.
In an almost Kenmotsu manifold we define a operator $h$ by $ h = \frac{1}{2}\pounds_{\xi}\varphi $ on $M$, where $\pounds_\xi$ is the Lie differentiation with respect to $\xi$. On an almost Kenmotsu manifold (see \cite{WL2, DP2})
\begin{eqnarray}\label{2.3}
h\xi  = 0, \hskip 0.3cm Tr~h = Tr~(h\varphi) =0, \hskip 0.3cm h\varphi = - \varphi h,
\end{eqnarray}
\begin{eqnarray}\label{2.4}
\nabla_{X}\xi = X - \eta(X)\xi - \varphi h X,
\end{eqnarray}
for any vector field $X$ on $M$; where $\nabla$ is the operator of covariant differentiation of $g$.\\

An almost Kenmotsu manifold $M^{2n+1}(\varphi,\xi,\eta,g)$ is said to be a generalized $(\kappa,\mu)$-almost Kenmotsu manifold if  $\xi$ belongs to the generalized $(\kappa,\mu)$-nullity distribution, i.e.,
\begin{eqnarray}\label{2.10}
R(X,Y)\xi = \kappa\{\eta(Y)X - \eta(X)Y\} + \mu\{\eta(Y)hX - \eta(X)hY\},
\end{eqnarray}
for all vector fields $X$, $Y$ on $M$, where $\kappa$, $\mu$ are smooth functions on  $M$.
An almost Kenmotsu manifold $M^{2n+1}(\varphi,\xi,\eta,g)$ is said to be a generalized $(\kappa,\mu)'$-almost Kenmotsu manifold if  $\xi$ belongs to the generalized $(\kappa,\mu)$-nullity distribution, i.e.,
\begin{eqnarray}\label{2.7}
R(X,Y)\xi = \kappa\{\eta(Y)X - \eta(X)Y\} + \mu\{\eta(Y)h' X - \eta(X)h' Y\},
\end{eqnarray}
for all vector fields $X$, $Y$ on $M$, where $h'=ho\varphi$ and $\kappa$, $\mu$ are smooth functions on $M$. Moreover, if both $\kappa$ and $\mu$ are constants in Eq$.$ (\ref{2.7}), then $M$ is called a $(\kappa,\mu)'$-almost Kenmotsu manifold.
Classifications of almost Kenmotsu manifolds with $\xi$ belong to $(\kappa$, $\mu)$-nullity distribution and $(\kappa,\mu)'$-nullity distribution have done by several authors. For more details, we refer the reader to \cite{WL2,DP2,DP3}. The following formulas are valid on a generalized $(\kappa,\mu)$ or $(\kappa,\mu)'$-almost Kenmotsu manifold (e.g., \cite{DP2})
\begin{eqnarray}\label{2.8}
h'^2 = (\kappa+1)\varphi^2~~~~~or,~~equivalently~~~~~ h^2 = (\kappa+1)\varphi^2.
\end{eqnarray}
\begin{eqnarray}\label{3.3}
Q\xi=2n\kappa \xi
\end{eqnarray}
Consider $X \in \mathcal{D}$ be an eigenvector of $h'$ with eigenvalue $\sigma$, where $\mathcal{D}$ is the distribution such that $\mathcal{D}=ker(\eta)$. It follows from (\ref{2.8}) that
 $\sigma ^2=-(\kappa+1)$ and therefore $\kappa\leq-1$ and $\sigma=\pm\sqrt{-\kappa-1}$. The equality holds if and only if $h=0$ (equivalently, $h'=0$). Thus, $h'\neq0$ if and only if $\kappa<-1$.

\section{Main Results}
Before entering into our main results we now recall the following

\begin{lem}\label{lem3.1}
(Proposition $3.2$ of \cite{DP3}) Let $M^{2n+1}(\varphi,\xi,\eta,g)$ be a generalized $(\kappa,\mu)'$-almost Kenmotsu manifold with $h\neq0$. Then
\begin{eqnarray}\label{A}
\xi(\lambda)=-\lambda (\mu+2), ~~~~~~~\xi(\kappa)=-2(\kappa+1)(\mu+2).
\end{eqnarray}
\end{lem}

Recently, Wang-Liu\cite{WL2} obtained the expression of Ricci operator on generalized $(\kappa,\mu)$ or $(\kappa,\mu)'$-almost Kenmotsu manifold.

\begin{lem}\label{lem3.2}
(Lemma $3.4$ of \cite{WL2}) Let $M^{2n+1}(\varphi,\xi,\eta,g)$ be a generalized $(\kappa,\mu)$-almost Kenmotsu manifold with $h'\neq0$. For $n>1$, the Ricci operator $Q$ of $M^{2n+1}$ can be expressed as
\begin{eqnarray}\label{C}
QX = -2nX + 2n(\kappa+1)\eta(X)\xi -2(n-1)h'X + \mu hX,
\end{eqnarray}
for any vector field $X$ on $M$. Also, the scalar curvature of $M$ is $2n(\kappa-2n)$.
\end{lem}
\begin{lem}\label{lem3.3}
(Lemma $3.3$ of \cite{WL2}) Let $M^{2n+1}(\varphi,\xi,\eta,g)$ be a generalized $(\kappa,\mu)'$-almost Kenmotsu manifold with $h'\neq0$. For $n>1$, the Ricci operator $Q$ of $M$ can be expressed as
\begin{eqnarray}\label{B}
QX = -2nX + 2n(\kappa+1)\eta(X)\xi -[\mu-2(n-1)h']X,
\end{eqnarray}
for any vector field $X$ on $M$. Further, if $\kappa$ and $\mu$ are constants and $n\geq1$, then $\mu=-2$ and hence
\begin{eqnarray}\label{2.9}
QX = -2nX + 2n(\kappa+1)\eta(X)\xi - 2nh'X,
\end{eqnarray}
for any vector field $X$ on $M$. In both case, the scalar curvature of $M$ is $2n(\kappa-2n)$.
\end{lem}

We now deduce the expression of the curvature tensor that satisfies the Miao-Tam equation

\begin{lem}\label{lem3.4}
Let a Riemannian manifold $(M^n, g)$ satisfies the Miao-Tam equation. Then the curvature tensor $R$ can be expressed as
\begin{eqnarray}\label{3.2}
&R(X,Y)D\lambda = (X\lambda)QY - (Y\lambda)QX + \lambda\{(\nabla_{X}Q)Y -(\nabla_{Y}Q)X \}\nonumber\\
& + (Xf)Y - (Yf)X,
\end{eqnarray}
for any vector fields $X$, $Y$ on $M$ and $f=-\frac{r\lambda + 1}{n-1}$.
\end{lem}
\textbf{Proof: }The equation (\ref{1.1}) can be exhibited as
\begin{eqnarray}\label{3.3A}
\nabla_{X}D\lambda = \lambda QX + (1+\triangle_{g}\lambda)X,
\end{eqnarray}
for any vector field $X$ on $M$. Now, tracing (\ref{1.1}) we obtain $\triangle_{g}\lambda = -\frac{r\lambda + n}{n-1}$. Then the Eq$.$ (\ref{3.3A}) transform into
\begin{eqnarray}\label{3.4}
\nabla_{X}D\lambda = \lambda QX + fX,
\end{eqnarray}
for any vector field $X$ on $M$.
Taking the covariant derivative of (\ref{3.4}) along an arbitrary vector field $Y$ on $M$, we obtain
\begin{align*}
\nabla_{Y}(\nabla_{X}D\lambda) = (Y\lambda)QX  + \lambda\{(\nabla_{Y}Q)X + Q(\nabla_{Y}X)\} + (Yf)X + f\nabla_{Y}X,
\end{align*}
for any vector field $X$ on $M$. Applying the preceding equation and (\ref{3.4}) in the well known expression of the curvature tensor
$R(X,Y) = [\nabla_{X},\nabla_{Y}]-\nabla_{[X,Y]},$
we obtain the required result.~~~~~~~~~~~~~~~~~~~~~~~~~~~~~~~~~~~~~~~~~~~~~~~~~~~~~~~~~~~~~~~~~~~~~~~~~~~$\square$\\

Dileo and Pastore (Proposition $4.1$ of \cite{DP2}) proved that on a $(\kappa,\mu)'$-almost Kenmotsu manifold, $\mu=-2$. So, we consider the Miao-Tam equation on $(\kappa,-2)'$-almost Kenmotsu manifold with $h'\neq0$ and prove

\begin{thm}
Let $M^{2n+1}(\varphi, \xi, \eta, g)$ be a $(\kappa,-2)'$-almost Kenmotsu manifold with $h'\neq0$. If there is a non-constant function $\lambda$ on $M$ satisfying the Miao-Tam equation, then $M^3$ is locally isometric to the Riemannian product $\mathbb{H}^2(-4)\times \mathbb{R}$, and for $n>1$, $M^{2n+1}$ is locally isometric to the warped products $\mathbb{H}^{n+1}(\alpha)\times_{f} \mathbb{R}^{n}$ or, $B^{n+1}(\alpha')\times_{f'}\mathbb{R}^n$; where $\mathbb{H}^{n+1}(\alpha)$ is the hyperbolic space of constant curvature $\alpha=-1-\frac{2}{n}-\frac{1}{n^2}$, $B^{n+1}(\alpha')$ ia a space of constant curvature $\alpha'=-1+\frac{2}{n}-\frac{1}{n^2}$, $f=ce^{(1-\frac{1}{n})t}$ and $f'=c'e^{(1+\frac{1}{n})t}$, with $c$, $c'$ positive constants.
\end{thm}
\textbf{Proof: }Since the scalar curvature of $(\kappa,-2)'-$almost Kenmotsu manifold with $h'\neq0$ is $2n(\kappa-2n)$ (follows from lemma \ref{lem3.3}), the equation \eqref{3.2} can be written as
\begin{eqnarray}\label{3.4A}
&R(X,Y)D\lambda = (X\lambda)QY - (Y\lambda)QX + \lambda\{(\nabla_{X}Q)Y -(\nabla_{Y}Q)X \}\nonumber\\
& + (2n-\kappa)\{(X\lambda)Y - (Y\lambda)X\},
\end{eqnarray}
for any vector fields $X$, $Y$ on $M$. Now, substituting $X$ by $\xi$ in \eqref{3.4A}, then taking its inner product with $\xi$ and using \eqref{3.3} we get
\begin{eqnarray}\label{E}
&g(R(\xi,Y)D\lambda,\xi) = \{2n(\kappa+1)-\kappa\}\{(\xi \lambda)\eta(Y) - (Y\lambda)\} \nonumber\\
&+ \lambda\{g((\nabla_{\xi}Q)Y,\xi) - g((\nabla_{Y}Q)\xi,\xi)\},
\end{eqnarray}
for all vector field $Y$ on $M$. Taking covariant derivative of (\ref{3.3}) along an arbitrary vector field $Y$ on $M$ we have
$(\nabla_{Y}Q)\xi + Q(\nabla_{Y}\xi)=2n\kappa \nabla_{Y}\xi.$
By virtue of (\ref{2.4}), the last equation gives
\begin{eqnarray}\label{D}
(\nabla_{Y}Q)\xi = 2n\kappa(Y-\varphi hY)-Q(X-\varphi hY),
\end{eqnarray}
for any vector field $Y$ on $M$.
Now, using (\ref{D}) and (\ref{2.3}) in (\ref{E}) provides
\begin{eqnarray}\label{F}
g(R(\xi,Y)D\lambda,\xi) &=& \{2n(\kappa+1)-\kappa\}\{(\xi\lambda)\eta(Y) - (Y\lambda)\},
\end{eqnarray}
for all vector field $Y$ on $M$. Further, taking the scalar product of (\ref{2.7}) with $D\lambda$ and then replacing $X$ by $\xi$ in the resulting equation and recalling $h' \xi=0$ and $g(Y, D\lambda)=Y\lambda$ gives
\begin{eqnarray*}
g(R(\xi,Y)D\lambda,\xi) = \kappa g(D\lambda-(\xi\lambda)\xi,Y) + 2 g(D\lambda,h' Y),
\end{eqnarray*}
for any vector field $Y$ on $M$, where we use $\mu=-2$. Combining the last two equations it follows that
\begin{eqnarray}\label{G}
n(\kappa+1)\{D\lambda-(\xi\lambda)\xi\} - h\varphi D\lambda = 0.
\end{eqnarray}
Now, operating the foregoing equation by $ h\varphi$ and using $h \xi=0$ yields
$n(\kappa+1) h\varphi D\lambda +  h^2 \varphi^2 D\lambda = 0.$
By virtue of (\ref{G}), (\ref{2.8}), $\varphi \xi=0$ and the $1$st Eq$.$ of (\ref{2.1}), the preceding Eq$.$ provides
$$n^2(\kappa+1)^2\{D\lambda-(\xi\lambda)\xi\} - (\kappa+1) \varphi^2D\lambda = 0.$$
Moreover, making use of (\ref{2.1}) the last equation reduces to
$(\kappa+1)\{n^2(\kappa+1)+1\}\{D\lambda-(\xi\lambda)\xi\} = 0.$
Since  $\kappa<-1$, the foregoing Eq$.$ gives
\begin{eqnarray}\label{3.5}
\{n^2(\kappa+1)+1\}\{D\lambda-(\xi\lambda)\xi\} = 0.
\end{eqnarray}
Thus, we have either $n^2(\kappa+1)+1=0$, or $n^2(\kappa+1)+1\neq0$.\\

\textbf{Case I: }In this case, we have $\kappa=-1-\frac{1}{n^{2}}$. For $n=1$, $\kappa=\mu=-2 $ and therefore from Theorem $4.2$ of Dileo and Pastore \cite{DP2} we deduce that $M^3$ is locally isometric to the Riemannian product $\mathbb{H}^2(-4)\times \mathbb{R}$ and for $n>1$, $M^{2n+1}$ is locally isometric to the warped products $\mathbb{H}^{n+1}(\alpha)\times_{f} \mathbb{R}^{n}$ or, $B^{n+1}(\alpha')\times_{f'}\mathbb{R}^n$; where $\mathbb{H}^{n+1}(\alpha)$ is the hyperbolic space of constant curvature $\alpha=-1-\frac{2}{n}-\frac{1}{n^2}$, $B^{n+1}(\alpha')$ is a space of constant curvature $\alpha'=-1+\frac{2}{n}-\frac{1}{n^2}$, $f=ce^{(1-\frac{1}{n})t}$ and $f'=c'e^{(1+\frac{1}{n})t}$, with $c$, $c'$ positive constants. \\

\textbf{Case II: }In this case, it follows from (\ref{3.5}) that $ D\lambda = (\xi \lambda)\xi $. Taking covariant derivative of $ D\lambda = (\xi \lambda)\xi $ along an arbitrary vector field $ X $ on $M$ and using (\ref{2.1}), (\ref{2.4}), we deduce
\begin{eqnarray}\label{3.5A}
\nabla_{X}D\lambda = X(\xi \lambda)\xi + (\xi \lambda)( X - \eta(X)\xi - \varphi h X).
\end{eqnarray}
Since the scalar curvature (from Lemma $3.3$) is $2n(\kappa-2n)$, the Eq$.$ (\ref{3.4}) can be written as
\begin{eqnarray}\label{3.5B}
&\nabla_{X}D\lambda = \lambda \{QX+(2n-\kappa)X\}-\frac{1}{2n}X,
\end{eqnarray}
for any vector field $X$ on $M$. Making use of (\ref{3.5A}) in (\ref{3.5B}) it follows that
\begin{eqnarray*}
&\lambda QX = \{(\kappa-2n)\lambda+(\xi\lambda)+\frac{1}{2n}\}X + X(\xi \lambda)\xi - (\xi \lambda)\{\eta(X)\xi + \varphi h X\},
\end{eqnarray*}
for any vector field $X$ on $M$. Comparing this with (\ref{2.9}) we deduce that
\begin{eqnarray}\label{3.6}
&\{\kappa\lambda+(\xi\lambda)+\frac{1}{2n}\}X+ X(\xi \lambda)\xi\nonumber\\
&+ \{(\xi\lambda)+2n\lambda\}h' X - \{2n(\kappa+1)\lambda+(\xi\lambda)\}\eta(X)\xi = 0,
\end{eqnarray}
for any vector field $X$ on $M$. Now, tracing (\ref{3.6}) over $X$ and noting that $Tr~h' = 0$, we have
\begin{eqnarray}\label{3.7}
&(2n+1)\{\kappa\lambda+(\xi\lambda)++\frac{1}{2n}\} + \xi(\xi\lambda) - \{2n(\kappa+1)\lambda+(\xi\lambda)\} = 0.
\end{eqnarray}
Next, substituting $X$ by $\xi$ in the equation (\ref{3.4}) and then taking its scalar product with $\xi$ yields
$\xi(\xi\lambda) = \lambda\{2n(\kappa+1)-\kappa\}-\frac{1}{2n}.$
By virtue of this, Eq$.$ (\ref{3.7}) takes the form
\begin{eqnarray}\label{3.8}
&\kappa\lambda+(\xi\lambda)+\frac{1}{2n} = 0.
\end{eqnarray}
Therefore, operating equation (\ref{3.6}) by $\varphi^2$ and taking into account (\ref{3.8}) we infer that $\{(\xi\lambda)+2n\lambda\}\varphi^2 h' X = 0$ for any vector field $X$ on $M^{2n+1}$. Further, making use of the first Eq$.$ of (\ref{2.1}), $h' = h\circ\varphi$, $h\varphi = -\varphi h$ and $\varphi\xi=0$, the last Eq$.$ reduces to
\begin{eqnarray}\label{3.9}
((\xi\lambda)+2n\lambda)h'X = 0,
\end{eqnarray}
for any vector field $X$ on $M$. By virtue of (\ref{3.8}), the equation (\ref{3.9}) gives $\{(2n-\kappa)\lambda-\frac{1}{2n}\} h'X = 0$ for any vector field $X$ on $M$.
Since $h'$ is non-vanishing, $\kappa<-1$. Thus, it follows that $\lambda = \frac{1}{2n(2n-\kappa)}$ and which is constant. This completes the proof.~~~~~~~~~~~~~~~~~~~~~~~~~~~~~~~~~~~~~~~~~~~~~~~~~~~~~~~~~~~~~~~~~~~~~~~~~~~$\square$\\

Next, we extend the last result for a generalized $(\kappa,\mu)'-$almost Kenmotsu manifold. Note that the metric satisfying the Miao-Tam equation has constant scalar curvature (see \cite{MT}). Further, the scalar curvature of generalized $(\kappa,\mu)'-$almost Kenmotsu manifold with $h'\neq0$ is $2n(\kappa-2n)$ (follows from lemma \ref{lem3.3}). Hence, it follows that $\kappa$ is constant. Therefore, from lemma \ref{lem3.1} we have $(\kappa+1)(\mu+2)=0$ . Since $\kappa<-1$, we must have $\mu=-2$. Thus from the last theorem we have

\begin{thm}
Let $M^{2n+1}(\varphi, \xi, \eta, g)$, $n>1$, be a generalized $(\kappa,\mu)'-$almost Kenmotsu manifold with $h'\neq0$. If there is a non-constant function $\lambda$ on $M$ satisfying the Miao-Tam equation, then $g$ is locally isometric to the warped products $\mathbb{H}^{n+1}(\alpha)\times_{f} \mathbb{R}^{n}$ or, $B^{n+1}(\alpha)\times_{f'}\mathbb{R}^n$; where $f=ce^{(1-\frac{1}{n})t}$ and $f'=c'e^{(1+\frac{1}{n})t}$, with $c$, $c'$ positive constants.
\end{thm}

Finally, we examine the existence of the solution of the Miao-Tam equation on a generalized $(\kappa,\mu)-$almost Kenmotsu manifold with $h\neq0$.

\begin{thm}
There does not exists any solution of the  Miao-Tam equation on a generalized $(\kappa,\mu)-$almost Kenmotsu manifold $M^{2n+1}(\varphi, \xi, \eta, g)$, $n>1$, with $h\neq0$.
\end{thm}
\textbf{Proof: }Suppose there exists a nontrivial smooth function $\lambda$ such that $(g,\lambda)$  is a solution of the Miao-Tam equation \eqref{1.1}.
Then it satisfies the curvature equation \eqref{3.2}.
Therefore, taking the scalar product of (\ref{3.2}) with $\xi$, and then using (\ref{3.3}), $r=2n(\kappa-2n)$ (follows from lemma \ref{lem3.2}) we achieve
\begin{eqnarray}\label{3.4E}
&g(R(X,Y)D\lambda,\xi) = \{2n(\kappa+1)-\kappa\}\{(X\lambda)\eta(Y) - (Y\lambda)\eta(X)\} \nonumber\\
&+ \lambda\{g(Y,(\nabla_{X}Q)\xi) - g(X,(\nabla_{Y}Q)\xi)\},
\end{eqnarray}
for all vector fields $X$, $Y$ on $M$. Since the metric satisfying (\ref{1.1}) has constant scalar curvature (see \cite{MT}) and $r=2n(\kappa-2n)$, it follows that $\kappa$ is constant. Hence the equation (\ref{D}) also valid here. Now, making use of (\ref{D}) and (\ref{2.3}) in (\ref{3.4E}), we immediately infer that
\begin{eqnarray}\label{3.E}
&g(R(X,Y)D\lambda,\xi) = \{2n(\kappa+1)-\kappa\}\{(X\lambda)\eta(Y) - (Y\lambda)\eta(X)\} \nonumber\\
&+\lambda\{g(Q\varphi hX,Y) - g(X,Q\varphi  hY)\},
\end{eqnarray}
for all vector fields $X$, $Y$ on $M$.
Next, putting $X=\varphi X$ and $Y=\varphi Y$ in (\ref{3.E}) and noting that $ g(R(\varphi X,\varphi Y)D\lambda,\xi) = 0$ (follows from (\ref{2.10})) and $h\varphi=-\varphi h$ we have
\begin{eqnarray*}
\lambda\{g(Qh\varphi^2X, \varphi Y)-g(\varphi X, Qh\varphi^2Y)\}=0,
\end{eqnarray*}
for all vector fields $X$, $Y$ on $M$. Therefore, using (\ref{C}) and $h\varphi=-\varphi h$ we have $\lambda\mu h^2 \varphi^2 X=0$ for any vector field $X$ on $M$. Thus, by virtue of (\ref{2.1}) and \eqref{2.8}, the last equation yields $(\kappa+1)\lambda\mu \varphi^2 X=0$ for any vector field $X$ on $M$. Since $\kappa<-1$, the foregoing eqiation provides $\lambda\mu=0$.\\

We suppose that $\lambda\neq0$ in some open set $\mathcal{O}$ of $M$. Then on $\mathcal{O}$, $\mu=0$.
Now, replacing $X$ by $\xi$ in (\ref{2.10}) and then taking the scalar product of the resulting Eq$.$ with $D\lambda$ gives
\begin{eqnarray}\label{M}
g(R(\xi,Y)D\lambda,\xi) = \kappa g(D\lambda-(\xi\lambda)\xi,Y).
\end{eqnarray}
Moreover, substituting $X$ by $\xi$ in (\ref{3.E}) and using \eqref{3.3}, $h\xi=0$ and $\varphi \xi=0$,  we have $g(R(\xi,Y)D\lambda,\xi) = \{2n(\kappa+1)-\kappa\}\{(\xi\lambda)\eta(Y) - (Y\lambda)\}$. Combining this with (\ref{M}) we obtain
\begin{eqnarray*}
(\kappa+1)\{D\lambda-(\xi\lambda)\xi\}=0.
\end{eqnarray*}
Since $h\neq0$, i.e., $\kappa<-1$, the last equation gives $D\lambda - (\xi \lambda)\xi=0$. Taking covariant derivative of $ D\lambda = (\xi \lambda)\xi $ along an arbitrary vector field $ X $ and using (\ref{2.1}), (\ref{2.4}), we deduce
\begin{eqnarray}\label{3.11}
\nabla_{X}D\lambda = X(\xi \lambda)\xi + (\xi \lambda)( X - \eta(X)\xi - \varphi h X).
\end{eqnarray}
Since the scalar curvature is $2n(\kappa-2n)$ (follows from lemma \ref{lem3.2}), the Eq$.$ (\ref{3.4}) transforms into
\begin{eqnarray}\label{3.5C}
&\nabla_{X}D\lambda = \lambda \{QX+(2n-\kappa)X\}-\frac{1}{2n}X.
\end{eqnarray}
Making use of (\ref{3.5C}) in (\ref{3.11}) it follows that
\begin{eqnarray*}
&\lambda QX = \{(\kappa-2n)\lambda+(\xi\lambda)+\frac{1}{2n}\}X + X(\xi \lambda)\xi - (\xi \lambda)(\eta(X)\xi + \varphi h X).
\end{eqnarray*}
By virtue of (\ref{C}), the last equation can written as
\begin{eqnarray}\label{3.12}
&\{\kappa\lambda+(\xi\lambda)+\frac{1}{2n}\}X+ X(\xi \lambda)\xi\nonumber\\
&+\{(\xi\lambda)+2(n-1)\lambda\}h\varphi X - \{2n(\kappa+1)\lambda+(\xi\lambda)\}\eta(X)\xi = 0.
\end{eqnarray}
Now, tracing (\ref{3.12}) over $X$ and noting that $Tr~h\varphi = 0$, we have
\begin{eqnarray}\label{3.19}
&(2n+1)\{\kappa\lambda+(\xi\lambda)+\frac{1}{2n}\} + \xi(\xi\lambda) - \{2n(\kappa+1)\lambda+(\xi\lambda)\} = 0.
\end{eqnarray}
Next, substituting $X$ by $\xi$ in (\ref{3.5C}) and using $Q\xi=2n\kappa \xi$ and then taking scalar product of the resulting equation with $\xi$, we get
$\xi(\xi\lambda) = \lambda\{2n(\kappa+1)-\kappa\}-\frac{1}{2n}.$
Making use of this, (\ref{3.19}) reduces to
\begin{eqnarray}\label{3.20}
\kappa\lambda+(\xi\lambda)+\frac{1}{2n} = 0.
\end{eqnarray}
Further, operating (\ref{3.12}) by $\varphi$ and using (\ref{3.20}) and  $h \varphi = -\varphi h$ we get $\{2(n-1)\lambda +(\xi\lambda)\} h\varphi^2 X = 0$. Moreover, using (\ref{2.1}) and recalling $\varphi\xi=0$, the last Eq$.$ provides
$\{2(n-1)\lambda +(\xi\lambda)\}hX = 0$.
Making use of (\ref{3.20}), the preceding Eq$.$ transform into $\{2(n-1)\lambda-\kappa\lambda-\frac{1}{2n}\} hX = 0$.
Since $h\neq0$, the last Eq$.$ shows that $2n\{2(n-1)-\kappa\}\lambda-1=0$. As $\kappa<-1$, this shows that $\lambda$ is constant. Therefore, it follows  from equation \eqref{1.1} that $S = - \frac{1}{\lambda} g$ on $\mathcal{O}$, since $\lambda\neq0$ on $\mathcal{O}$. The $g$-trace gives $- \frac{1}{\lambda}= \frac{r}{2n+1} = \frac{2n(\kappa-2n)}{2n+1}$. Therefore, $S=\frac{2n(\kappa-2n)}{2n+1} g$. Since $Q\xi = 2n\kappa\xi$, from the foregoing equation we deduce that $\kappa = -1$ on $\mathcal{O}$, which is a contradiction. Hence $\lambda$ is trivial on $M$. This completes the proof. ~~~~~~~~~~~~~~~~~~~~~~~~~~~~~~~~~~~~~~~~~~~~~~~~~~~~~~~~~~~~~~~~~~~~~~~~~~~$\square$
%

\noindent\\
\textbf{Acknowledgments:}
The author D. S. Patra is financially supported by the Council of Scientific and Industrial Research, India (grant no. 17-06/2012(i)EU-V).

$^1$
Department of Mathematics, \\
Jadavpur University,  \\
188. Raja S. C. Mullick Road,\\
Kolkata:700 032, INDIA \\
E-mail: dhritimath@gmail.com\\

$^2$
Department of Mathematics, \\
Chandernagore College\\
Hooghly: 712 136 (W.B.), INDIA\\
E-mail: aghosh\_70@yahoo.com\\

\end{document}